\documentclass{amsart}
\usepackage[utf8]{inputenc}
\usepackage{amsmath, amsthm}
\usepackage{amssymb}
\usepackage[textsize=footnotesize]{todonotes}
\usepackage{tikz-cd}
\usepackage{enumerate}
\usepackage[utf8]{inputenc}
\usepackage{amsmath}
\usepackage{newunicodechar}
\usepackage{comment}

\usepackage{url}
\usepackage{hyperref}

\DeclareFontFamily{U}{min}{}
\DeclareFontShape{U}{min}{m}{n}{<-> udmj30}{}

\newtheorem{theorem}{Theorem}
\newtheorem*{theorem*}{Theorem}
\newtheorem*{corollary*}{Corollary}

\newtheorem{lemma}[theorem]{Lemma}
\newtheorem{prop}[theorem]{Proposition}

\newtheorem{claim}[theorem]{Claim}
\newtheorem{corollary}[theorem]{Corollary}

\newcommand{\la}{\langle}
\newcommand{\ra}{\rangle}
\newcommand{\rest}{\restriction}

\newcommand{\p}{\mathcal{P}}

\newcommand{\bbP}{\mathbb{P}}

\newcommand{\calA}{\mathcal{A}}
\newcommand{\calB}{\mathcal{B}}

\DeclareMathOperator{\dom}{dom}

\DeclareMathOperator{\cf}{cf}

\DeclareMathOperator{\crit}{crit}
\DeclareMathOperator{\ot}{ot}

\DeclareMathOperator{\cof}{cof}

\author{Monroe Eskew}\address{Kurt G\"odel Research Center,
University of Vienna,
Vienna, Austria.}
\email{monroe.eskew@univie.ac.at}
\title{Forcing square with small conditions}
\date{}

\thanks{This research was supported by Austrian Science Fund (FWF) Project PIN1355423.}

\begin{document}
\maketitle

\begin{abstract}
    We show that if $\kappa^{<\kappa} = \kappa$, then there is a $\kappa$-strategically closed, $\kappa^+$-c.c.\ partial order that forces $\square_\kappa$.  We discuss consequences for generic large cardinals.
\end{abstract}

\section{Introduction}
The square principle $\square_\kappa$ was introduced by Jensen \cite{jensen}, who showed that in G\"odel's constructible universe $L$, $\square_\kappa$ holds for all infinite cardinals $\kappa$.  The principle asserts that there is a sequence $\la C_\alpha : \alpha \in \lim(\kappa^+)\ra$ such that:
\begin{enumerate}[(1)]
    \item For each limit ordinal $\alpha<\kappa^+$, $C_\alpha$ is a closed unbounded (club) subset of $\alpha$ of ordertype at most $\kappa$.
    \item If $\beta \in \lim(C_\alpha)$, then $C_\beta = C_\alpha\cap\beta$.
\end{enumerate}
Extending Jensen's result about $L$, $\forall\kappa\square_\kappa$ has been shown to hold in many canonical inner models \cite{SZ-square}.
Owing to the upwards-absoluteness of $\square_\kappa$ to any outer model with the same $\kappa^+$, square principles have been pivotal in inner model theory, with their failure implying the consistency of various large cardinals \cite{schimmerling}.  

Jensen also showed that $\square_\kappa$ can also be forced over any model of set theory by a $(\kappa+1)$-strategically closed forcing that adds the sequence by initial segements.  In particular, Jensen's forcing adds no $\kappa$-sequences of ordinals, is of size $2^\kappa$, and is not $2^\kappa$-c.c.

$\square_\kappa$ is a kind of incompactness property for $\kappa^+$, as a $\square_\kappa$-sequence
cannot be extended to include $\kappa^+$ itself, even with the ordertype requirement dropped at that point.  A club $C\subseteq\kappa^+$ will have a tail of limit points $\alpha$ such that $\ot(C\cap\alpha)>\kappa$, and thus for a $\square_\kappa$-sequence $\la C_\alpha : \alpha\in \lim(\kappa^+)\ra$, $C\cap\alpha\not= C_\alpha$ for most $\alpha$.

Solovay \cite{solovay} showed that $\square_\lambda$ fails for all $\lambda$ above a strongly compact cardinal.  This motivates the question of whether some analog of a large cardinal property that can hold at small cardinals may negate squares.  In particular, does some generic large cardinal property of $\kappa^+$ imply the sort of compactness for $\kappa^+$ that would rule out $\square_\kappa$?

It was shown in \cite{localsat} that, assuming the consistency of a huge cardinal, it is consistent with ZFC that for all $\kappa$, $\square_\kappa$ holds, while the nonstationary ideal on $\kappa^+$ is somewhere saturated.  On the other hand, some related generic large cardinal properties, such as the Chang Conjecture $(\kappa^{++},\kappa^+)\twoheadrightarrow(\kappa^+,\kappa)$, are known to imply the failure of $\square_{\kappa^+}$ \cite{foreman2002}.
The strategy of \cite{localsat} for getting $\square_\kappa$ along with a saturated ideal on $\kappa^+$ was to start with a model with a robustly saturated ideal on $\kappa^+$, force $\square_{\kappa^+}$ with Jensen's initial segments forcing while noting the preservation of the saturation property, and then use Jensen's forcing for $\square_\kappa$, arguing that the companion ``threading'' forcing, along with the enhanced strategic closure we get from $\square_{\kappa^+}$ holding in the background, enables the construction of a master filter.
This strategy obviously doesn't work to get $\square_\kappa$ along with principles that contradict $\square_{\kappa^+}$.

On the other hand, it was shown by Shelah and Stanley \cite{shelah-stanley} that under CH, $\square_{\omega_1}$ can be forced by a countably closed, $\omega_2$-c.c.\ forcing.  This allows one to easily show that many generic large cardinal properties of $\omega_2$, including  $(\omega_3,\omega_2)\twoheadrightarrow(\omega_2,\omega_1)$, are consistent with $\square_{\omega_1}$.  However, it was unclear how to generalize this to higher cardinals. The present work clears this up:

\begin{theorem}\label{mainthm}
    If $\kappa^{<\kappa} = \kappa$, then there is a $\kappa$-strategically closed, $\kappa^+$-c.c.\ forcing that introduces a $\square_\kappa$-sequence.  Moreover, the $\kappa^+$-c.c.\ of the forcing is upwards-absolute to any outer model in which $\kappa$ and $\kappa^+$ remain regular.
\end{theorem}

Besides the forcings of Jensen and Shelah-Stanely, there have been several other notions developed for forcing squares.  Dolinar-D\u{z}amonja \cite{dolinar}, Krueger \cite{krueger}, and Neeman \cite{neeman} each give forcings for $\square_{\omega_1}$ with \emph{finite} conditions.  The posets of Krueger and Neeman are moreover $\omega_2$-c.c.  Hayut, Rinot, and You \cite{hayut2026notionsforcingaddsquare} show that for every regular $\kappa$, the principle $\diamondsuit_{\kappa^+}(\cof(\kappa))$ implies the existence of a $\kappa^+$-c.c., $\kappa^+$-distributive poset forcing $\square_\kappa$; the forcing is a particular $\kappa^+$-Suslin tree.  Our forcing here is different in that it adds many subsets of $\kappa$ and only requires a bit of GCH to work.

\begin{corollary}\label{glccor}
    Suppose $\kappa^{<\kappa} = \kappa$.  Assume either:
    \begin{enumerate}[(a)]
        \item There is $\kappa^+$-complete, $\kappa^{++}$-saturated ideal over a set $Z$.
        \item There is a normal, $\kappa^+$-complete, $\lambda^+$-saturated ideal on $[\lambda]^{<\kappa^+}$.
        \item There is a normal, $\kappa^+$-complete, $\lambda$-saturated ideal on $[\lambda]^{\kappa^+}$.
        \item $\kappa^+\leq\lambda_0<\lambda_1$ and $(\lambda_1,\lambda_0)\twoheadrightarrow(\kappa^+,\kappa)$.
    \end{enumerate}
            Then there is a cardinal-preserving poset that forces $\square_\kappa$ while preserving the assumed statement.
\end{corollary}

\section{The forcing and its properties}

Let us define a forcing $\bbP$ with conditions $p = \la c^p_\alpha : \alpha \in d^p \ra$, where:
\begin{enumerate}[(C1)]
    \item $d^p \in [\kappa^+]^{<\kappa}$, and $d^p$ is closed set of limit ordinals. 
    \item For each $\alpha \in d^p$, $c^p_\alpha$ is a closed subset of $\alpha$ of size $<\kappa$.
    \item If $\alpha \in d^p$ and $\cf(\alpha)<\kappa$, then $\sup(c^p_\alpha) = \alpha$.
    \item If $\alpha \in d^p$ and $\beta \in \lim(c^p_\alpha)$, then $\beta \in d^p$ and $c^p_\beta = c^p_\alpha \cap \beta$.
\end{enumerate}
We put $q \leq p$ when:
\begin{enumerate}[(O1)]
    \item $d^q \supseteq d^p$.
    \item For all $\alpha \in d^p$, $c^q_\alpha$ is an end-extension of $c^p_\alpha$.
    \item If $\alpha \in d^p$ and $\cf(\alpha)=\kappa$, then $\min(c^q_\alpha \setminus c^p_\alpha) > \sup(d^p \cap \alpha)$.
\end{enumerate}

\begin{lemma}
    For any $p \in \bbP$, any $\alpha \in \lim(\kappa^+)$, and any $\gamma<\alpha$, there is $q \leq p$ such that $\alpha \in d^q$ and $\sup(c^q_\alpha)>\gamma$.
\end{lemma}

\begin{proof}
    Let $p,\alpha,\gamma$ be as hypothesized.  
    If $\cf(\alpha)<\kappa$ and $\alpha\notin d^p$, then by (C1), $d^p \cap \alpha$ is bounded in $\alpha$, and so $c^q_\alpha$ can be chosen to trivially satisfy coherence with the other $c^p_\beta$'s.  If $\cf(\alpha) = \kappa$, then we can find $\delta$ such that 
    $$\max\{\gamma,\sup(d^p \cap \alpha),\sup (c^p_\alpha) \} < \delta<\alpha$$
    and put $c^q_\alpha = c^p_\alpha \cup \{\delta\}$ if $c^p_\alpha$ is defined, and otherwise put $c^q_\alpha = \{\delta\}$.
\end{proof}

Note that the above lemma implies that a sufficiently generic filter for $\bbP$ will give a $\square_\kappa$-sequence, provided that adjoining it preserves $\kappa^+$.  We also want to show that $\bbP$ adds no ${<}\kappa$-sequences.  We say that a poset is $\kappa$-strategically closed when Player II has a winning strategy in the game where two players take turns playing a descending sequence of conditions, with Player I starting and Player II playing at all limit rounds.  Player II wins when the game continues for $\kappa$-many rounds.  For notational convenience, we will index the first round by $1$, so that Player II plays at all even rounds.

For a set of ordinals $X$ and an ordinal $\alpha$, we will say, ``$X$ has a gap at $\alpha$,'' when $\sup(X\cap\alpha)<\alpha$.  We will say, ``$Y$ narrows the gaps of $X$ on $Z$,'' if whenever $\alpha \in Z$ and $X$ has a gap at $\alpha$, then $\sup(Y\cap\alpha) > \sup(X\cap\alpha)$.

\begin{lemma}
    $\bbP$ is $\kappa$-strategically closed.
\end{lemma}

\begin{proof}
    We describe a strategy for Player II.
    At successor rounds where II plays, if the previous play by Player I was $p$, then II plays $q$ such that $d^q$ narrows the gaps of $d^p$ on limits of limits below $\kappa^+$.  This is possible since if $d^p$ has a gap at $\alpha$ which is limit of limit, then II can adjoin $$\la \sup(d^p\cap\alpha)+\omega, \{\sup(d^p\cap\alpha)+n : n < \omega \} \ra.$$ 
    Since $\alpha\mapsto\sup(d^p\cap\alpha)$ takes ${<}\kappa$-many values, adding these all simultaneously yields a condition.
    
    At limit stages, the strategy is as follows.
    Suppose $\lambda<\kappa$ is a limit and the moves so far have been $\la p_i : i < \lambda \ra$.
    Let $d_0 = \bigcup_i d^{p_i}$.  
    For each $\alpha \in d_0$, let $\delta_\alpha = \sup(\bigcup_i c^{p_i}_\alpha)$, which is equal to $\alpha$ if and only if $\cf(\alpha)<\kappa$.  Note that if $\cf(\alpha) = \kappa$ and $\delta_\alpha \not= \max(c^{p_i}_\alpha)$ for some $i<\lambda$,
    then $\delta_\alpha \notin d_0$, because by (O3), $\beta \in d^{p_i}\cap\alpha$ implies $\min(c^{p_{i+1}}_\alpha \setminus c^{p_{i}}_\alpha)>\beta$.  Let 
    $$e = \{ \delta_\alpha :  \cf(\alpha) = \kappa \text{ and }\la c^{p_i}_\alpha : i < \lambda \ra \text{ is not eventually constant}\}.$$
    Note that if $\delta_\alpha = \delta_\beta \in e$, then $\alpha = \beta$.  This is because if $\alpha<\beta$ are in $d_0$, then $\sup(c^{p_i}_\beta) > \alpha$ for some $i<\lambda$.
    Let $d_1 = d_0 \cup e$, and let $d_2$ be the closure of $d_1$.

    If $\alpha \in d_0$, define $c^{*}_\alpha$ as the closure in $\alpha$ of $\bigcup_i c^{p_i}_\alpha$.  If $\alpha \in e$, define $c^{*}_\alpha = \bigcup_i c^{p_i}_\beta$, where $\beta$ is the unique ordinal such that $\delta_\beta = \alpha$.

    \begin{claim}\label{stable}
        Suppose $\alpha \in d_2 \setminus d_1$.  Then for all $\beta \in d_0$, $\alpha\notin\lim(c^{*}_\beta)$. Moreover, if $\beta>\alpha$, then there is $i<\lambda$ such that $c^{*}_\beta\cap\alpha = c^{p_i}_\beta\cap\alpha$.
    \end{claim}
    \begin{proof}
        If $\cf(\beta)<\kappa$, then there is $i<\lambda$ such that $c^{*}_\beta = c^{p_j}_\beta$ for $i\leq j <\lambda$.  If $\alpha\in\lim(c^{p_i}_\beta)$, then we would have $\alpha \in d^{p_i}$, contrary to assumption.

        If $\beta>\alpha$, $\cf(\beta) = \kappa$, and $c^{p_i}_\beta\cap\alpha$ is not eventually constant, then by (O3), $\{\max(c^{p_i}_\beta\cap\alpha):i<\lambda\}$ must be cofinal in $\alpha$, and $\alpha = \delta_\beta \in e$, contrary to assumption.  Thus there is $i<\lambda$ such that $c^{*}_\beta\cap\alpha = c^{p_j}_\beta\cap\alpha$ for $i\leq j<\lambda$.  Since $\alpha\notin d_0$,
        $\max(c^{*}_\beta\cap\alpha)<\alpha$.
    \end{proof}

    For $\alpha \in d_2 \setminus d_1$, let $\xi_\alpha$ be the least ordinal in $d_0 \setminus \alpha$ of cofinality $\kappa$, if such exists.  If $\xi_\alpha$ is defined, let $i_0^\alpha$ be the least $i<\lambda$ such that $\xi_\alpha \in d^{p_i}$ and $c^{p_i}_{\xi_\alpha}$ already possesses the largest limit point $\eta_\alpha$ of $c^*_{\xi_\alpha}\cap\alpha$, if such exists.  If $i^\alpha_0$ or $\eta_\alpha$ are undefined, set these variables equal 0.  Note that for $ i^\alpha_0<j<\lambda$, $\max(d^{p_j}\cap\alpha)>\eta_\alpha$.
    Then let $ c^{*}_\alpha = \{ \max(d^{p_j}\cap\alpha) : i_0^\alpha < j <\lambda$ and $j$ is even$\}$.

    We claim that $\{ \la \alpha,c^{*}_\alpha \ra : \alpha \in d_2\}$ is a condition below $p_i$ for all $i<\lambda$.  
    (C1) is clear.  
    (C2) is obvious for $\alpha\in d_1$.  
    (C3) holds for $\alpha\in d_0$ automatically and for $\alpha \in d_2 \setminus d_0$ by the definitions of $c^*_\alpha$ for those points in $e$ and those in $d_2\setminus d_1$.
    (C4) holds automatically for $\alpha\in d_0$ such that $\cf(\alpha)<\kappa$.  It holds for $\alpha \in d_0$ such that $\cf(\alpha) = \kappa$ by the addition of $e$ and the sets assigned at those points.  For $\alpha \in e$, it holds because if $\alpha = \delta_\beta$ and $\gamma \in \lim(c^*_\alpha)$, then for some $i<\lambda$, $\gamma \in \lim(c^{p_i}_\beta)$, and thus $\gamma \in d^{p_i}\subseteq d_0$ and $c^*_\alpha\cap\gamma = c^{p_i}_\gamma$.

    Thus it remains to check (C2) and (C4) for $\alpha \in d_2 \setminus d_1$.  First assume $\xi_\alpha$ is undefined.  Assuming the same strategy has been followed at previous stages, if $\gamma\in\lim(c^*_\alpha)$, then at the limit stage $i = \ot(c^*_\alpha\cap\gamma)$, $\gamma$ cannot be in $d^{p_j}$ for some $j<i$ because of the narrowing of gaps at successor rounds.  Thus Player II added $c^{p_i}_\gamma = \{ \max( d^{p_j}\cap\gamma) : j < i$ is even$\} = c^*_\alpha\cap\gamma$, verifying (C4) for $\alpha$.  Since $\gamma = \sup(\bigcup_{j<i} d^{p_j} \cap \alpha)$,
    Player II did not add any points to $d^{p_i} \cap (\gamma,\alpha)$, and so $\gamma = \max(d^{p_i}\cap\alpha) \in c^*_\alpha$, showing (C2) for $\alpha$.
    
    Next assume $\xi_\alpha>\alpha$ is defined.  Then for limit $j$ such that $i_0^\alpha<j<\lambda$, $\xi_\alpha$ is also the least ordinal of cofinality $\kappa$ in $\bigcup_{k<j} d^{p_k}$ above $\gamma = \sup_{k<j}(\max(d^{p_k}\cap\alpha))$.
    Since $\eta_\alpha<\gamma$, $i^\alpha_0$ is the least $i < j$ such that $c^{p_i}_{\xi_\alpha}$ possesses $\eta_\alpha$, which is also the largest limit ordinal of $\bigcup_{k<j} c^{p_k}_{\xi_\alpha}\cap\gamma$.  Note that by (O3) and the minimality of $\xi_\alpha$, there is no $\beta>\gamma$ such that $\gamma = \sup(\bigcup_{k<j} c^{p_k}_\beta\cap\gamma)$. 
    Thus assuming II has followed the same strategy earlier, they added $c^{p_j}_\gamma = \{ \max( d^{p_k}\cap\gamma) : i_0^\alpha<k < j$ and $k$ is even$\}= c^*_\alpha\cap\gamma$ at stage $j$, verifying (C4).  Since $\gamma = \sup(\bigcup_{k<j} d^{p_j} \cap \alpha)$, and no points in $(\gamma,\alpha)$ were limit points of $\bigcup_{k<j} c^{p_k}_\beta$ for some $\beta>\alpha$, Player II did not add any points to $d^{p_j} \cap (\gamma,\alpha)$, and so $\gamma = \max(d^{p_j}\cap\alpha) \in c^*_\alpha$, showing (C2).

    The ordering criterion (O1) is clear.  (O2) holds because for $\alpha \in d_0$, if $\cf(\alpha)<\kappa$ then $c^*_\alpha = c^{p_i}_\alpha$ for some $i<\lambda$, and if $\cf(\alpha) = \kappa$, then $c^*_\alpha$ is the union of a chain of end-extensions.  For (O3), note that if $\alpha \in d^{p_i}\cap\cof(\kappa)$, then $\min(c^*_\alpha\setminus c^{p_i}_\alpha) =\min(c^{p_{i+1}}_\alpha\setminus c^{p_i}_\alpha) > \sup(d^{p_i}\cap\alpha)$.
    \end{proof}

    Following the terminology of \emph{layered ideals} appearing in \cite{fms2}, the author chose in \cite{localsat} to say that a poset $\bbP$ is \emph{$S$-layered} for a set $S$ contained in a regular cardinal $\lambda$ when there is a club $C \subseteq \lambda$ and a decomposition of $\bbP$ as a $\subseteq$-increasing union of regular suborders $\la \bbP_\alpha : \alpha < \lambda \ra$, each of size $<\lambda$, that is continuous at all points in $S\cap C$.  It is easy to show that if $\lambda$ is regular, $S\subseteq\lambda$ is stationary, and $\bbP$ is $S$-layered, then $\bbP$ is $\lambda$-c.c.  Note that the notion of being $S$-layered is upwards-absolute, and thus the relevant chain condition is preserved in outer models where $\lambda$ remains regular and $S$ remains stationary.

    For $\delta<\kappa^+$, let $\bbP_\delta = \{ p \in \bbP : d^p \subseteq\delta \}$.  Note that if $\cf(\delta) = \kappa$, then $\bbP_\delta = \bigcup_{\alpha<\delta} \bbP_\alpha$.  For $p \in \bbP$ and $\delta<\kappa^+$, let $p \rest \delta$ denote $\la c^p_\alpha : \alpha \in d^p \cap\delta\ra$.  Note that if $\cf(\delta)=\kappa$, then for all $p \in \bbP$, $d^p\cap\delta$ is closed in $\kappa^+$, and thus $p \rest \delta \in \bbP$.
    
    \begin{lemma}
    If $\delta<\kappa^+$ and $\cf(\delta)=\kappa$, then $p \mapsto p \rest\delta$ is a projection from $\bbP$ to $\bbP_\delta$, and $\bbP_\delta$ is a regular suborder of $\bbP$.
    Thus if $\kappa^{<\kappa} = \kappa$, then $\bbP$ is $(\kappa^+ \cap \cof(\kappa))$-layered.
    \end{lemma}
    \begin{proof}
    It is easy to check that for $p,q\in\bbP$, $p \leq p \rest\delta$, and if $q \leq p$, then $q \rest \delta \leq p \rest \delta$. 
    Assume $p \in \bbP$ and $q \leq p\rest \delta$ is in $\bbP_\delta$.  We claim that $q \cup \{ \la\alpha,c^p_\alpha\ra : \alpha \in d^p \setminus \delta \}$ is a condition below $p$.  The only things that need checking are (C4) and (O3).  For (C4), notice that if $\alpha \in d^p \setminus \delta$ and $\beta \in \lim(c^p_\alpha)\cap \delta$, then $\cf(\beta)<\kappa$, $\beta \in d^p$, and $c^p_\beta = c^p_\alpha \cap\beta$.  Thus $c^q_\beta = c^p_\beta$.  
    For (O3), note that if $c^q_\beta \not= c^p_\beta$, then $\beta < \delta$ and $\min(c^q_\beta \setminus c^p_\beta) > d^{p \rest \delta}\cap\beta = d^p \cap\beta$.
    
    Therefore, if $\calA\subseteq\bbP_\delta$ is a maximal antichain of $\bbP_\delta$, then for any $p \in \bbP$, there is $q\in \bbP_\delta$ below $p\rest \delta$ and some $a \in \calA$, and thus there is some $p'\leq p$ such that $p' \rest \delta = q \leq a$.
\end{proof}

This completes the proof of Theorem~\ref{mainthm}, noting that the expression ``$\kappa^+\cap\cof(\kappa)$'' defines the same set in any outer model in which $\kappa$ and $\kappa^+$ remain regular.

We note the following about the quotient forcings $\bbP/\bbP_\delta$:

\begin{prop}
    Suppose $\delta<\kappa^+$, $\cf(\delta) = \kappa$, and $G_\delta\subseteq\bbP_\delta$ is generic over $V$. Then in $V[G_\delta]$, the quotient forcing $\bbP/G_\delta$ is not $\sigma$-strategically closed.
\end{prop}

\begin{proof}
  Let $\dot C_\delta^G$ be a $\bbP$-name for $\bigcup \{ c^p_\delta : p \in \dot G\}$.  First we claim that $\Vdash_\bbP \dot C_\delta^G\notin V[\dot G_\delta]$.  Suppose otherwise, and let $p,\tau$ be such that $\tau$ is a $\bbP_\delta$-name and $p \Vdash_\bbP \tau = \dot C^G_\delta$.  We may assume that $\delta \in \dom p$.  There are $p_0,p_1 \leq p$ that are the same as $p$ on $d^p \setminus \{\delta\}$ and such that for two ordinals $\alpha<\beta<\delta$, $c^{p_0}_\delta = c^p_\delta \cup \{\alpha \}$ and $c^{p_1}_\delta = c^p_\delta \cup \{\beta \}$.  Hence $p_0 \Vdash \alpha \in \tau$ and $p_1 \Vdash \alpha\notin\tau$.   But the statement ``$\alpha\in\tau$'' can be written in the forcing language of $\bbP_\delta$, and for any such atomic sentence $\varphi$, if $q \Vdash_\bbP \varphi$, then also $q\rest \delta \Vdash_{\bbP_\delta} \varphi$.
Hence $p_0 \rest \delta$ and $p_1 \rest \delta$ must make the same decision on whether $\alpha\in\tau$ as do $p_0$ and $p_1$ respectively.  But $p_0 \rest \delta = p_1 \rest \delta$, a contradiction.

Suppose then that $\bbP/G_\delta$ is $\sigma$-strategically closed.  Then $\kappa$ must be uncountable.  Let $p \in \bbP/G_\delta$, and using the fact that $\Vdash_\bbP \dot C^G_\delta \notin V[G_\delta]$, find $p_0,p_1 \leq p$ in $\bbP/G_\delta$ and $\alpha<\delta$ such that $p_0 \Vdash \alpha \in\dot C^G_\delta$, and $p_1 \Vdash \alpha \notin \dot C^G_\delta$.  Using a winning strategy $\Sigma$ for Player II, build a sequences $\la\la q_0^n,q_1^n\ra :  n < \omega \ra$ and $\la\la\beta_0^n,\beta_1^n :  n <\omega \ra$ such that:
\begin{enumerate}
    \item For $i<2$ and $n<\omega$, $q_i^0 = p_i$, and $\la q^m_i : 1\leq m < n \ra$ is a run of the game following $\Sigma$.
    \item For $i<2$ and odd $n$, $q_i^n \Vdash \beta_i^n \in \dot C^G_\delta$.
    \item For all $n<\omega$, $\beta_0^n<\beta_1^n < \beta_0^{n+1}$.
\end{enumerate}
Let $\gamma = \sup \{ \beta^n_0 : n < \omega \}$.  For $i<2$, since $\la q^n_i : n < \omega \ra$ follows $\Sigma$, there is a lower bound $q^\omega_i$, which forces that $\gamma$ is a limit point of $\dot C^G_\delta$.  Hence $q^\omega_i$ forces that $\dot C^G_\delta\cap \gamma = c^r_\gamma$ for some $r \in G_\delta$.
However, all $r \in G_\delta$ with $\gamma \in \dom r$ have the same value for $c^r_\gamma$.  This is a contradiction, as $q^\omega_0$ and $q^\omega_1$ disagree about whether $\alpha\in\dot C^G_\delta\cap \gamma$.
\end{proof}

\section{Consequences for generic large cardinals}

In this section we discuss the proof of Corollary~\ref{glccor}.  For a poset $\bbP$, let Let $\calB(\bbP)$ denote its Boolean completion.  The main tool is the Duality Theorem of \cite{foreman-duality}, which has the following as a special case:
\begin{theorem*}[Foreman]
\label{duality}
    Suppose $I$ is a $\kappa$-complete precipitous ideal on a set $Z$, and $\bbP$ is a $\kappa$-c.c.\ forcing.  Let $\bar I$ denote the ideal generated by $I$ in $V^{\bbP}$, and let $j$ denote the generic ultrapower embedding given by forcing with $\p(Z)/I$.    Then:
    $$\calB ( \bbP * \p(Z)/\bar I ) \cong \calB( \p(Z)/I * j(\bbP))$$
\end{theorem*}

Let us prove the first claim of Corollary~\ref{glccor}.  Suppose $\kappa^{<\kappa} = \kappa$ and $I$ is a $\kappa^+$-complete, $\kappa^{++}$-saturated ideal over a set $Z$.  Let $\bbP$ be the forcing from the previous section.  Let $G \subseteq \p(Z)/I$ be generic and let $j : V \to M \subseteq V[G]$ be the generic ultrapower embedding.  Then $\crit(j) = \kappa^+$ and $j(\kappa^+) = (\kappa^{++})^V$.  Further, $M \models j(\bbP)$ is $(j(\kappa^+)\cap\cof(\kappa))$-layered.  By the $\kappa$-closure of $M$ (see \cite{foremanhandbook}), ``$\kappa^+ \cap \cof(\kappa)$'' denotes the same set whether computed in $M$ or $V[G]$.  Hence $j(\bbP)$ has the $(\kappa^{++})^V$-c.c.\ in $V[G]$, and thus $\p(Z)/I * j(\bbP)$ is $\kappa^{++}$-c.c.\ in $V$.  By Foreman's Theorem, $\bbP$ forces that $\bar I$ is $\kappa^{++}$-saturated.

For the second and third claims of Corollary~\ref{glccor}, the argument is nearly the same as above.  The hypothesis gives that $j(\kappa^+) = \lambda^+$ in the case of the second claim and that $j(\kappa) = \lambda$ in the third.  $j(\bbP)$ will have the desired layering property in $M$ and $V[G]$ for the same reasons.
For the final claim, we apply the following folklore fact, a proof of which can be found in \cite{gcc}.
\begin{lemma}
    If $(\lambda_1,\lambda_0)\twoheadrightarrow(\kappa^+,\kappa)$, then this is preserved by any $\kappa^+$-c.c.\ forcing.
\end{lemma}

\bibliographystyle{amsplain.bst}
\bibliography{references.bib}

    \end{document}